\documentclass[12pt]{amsart}
\usepackage[a4paper,left=2cm,right=1.2cm,top=3cm,bottom=4cm]{geometry}
\usepackage{amsmath,amstext, amsthm,
amsbsy,amssymb,marvosym,fancyhdr,graphicx,amscd,amsfonts,latexsym,delarray,stackrel,lineno,color,cite,appendix,xcolor,url,hyperref}
\usepackage{wasysym,multirow,tabularx}
\usepackage{subcaption}

\usepackage{srcltx}
\usepackage{mathrsfs}
\usepackage{float}

\usepackage[all]{xy}
\usepackage{t1enc}
\usepackage{mathrsfs}
\usepackage{pifont}

\usepackage[table]{xcolor}

\usepackage{mathpple}

\usepackage[T1]{fontenc}

\definecolor{Green}{rgb}{0,1,0}
\definecolor{Blue}{RGB}{0,0,191}
\definecolor{mathmodecolor}{RGB}{0,102,0}
\definecolor{keywordcolor}{RGB}{0,51,151}
\definecolor{sourcebackgroundcolor}{RGB}{255,247,223}
\definecolor{unixagred}{RGB}{255,0,0}
\definecolor{lightgray}{RGB}{191,191,191}
\definecolor{green}{RGB}{1,191,191}

\newcommand*\patchAmsMathEnvironmentForLineno[1]{%
  \expandafter\let\csname old#1\expandafter\endcsname\csname #1\endcsname
  \expandafter\let\csname oldend#1\expandafter\endcsname\csname end#1\endcsname
  \renewenvironment{#1}%
     {\linenomath\csname old#1\endcsname}%
     {\csname oldend#1\endcsname\endlinenomath}}%
\newcommand*\patchBothAmsMathEnvironmentsForLineno[1]{%
  \patchAmsMathEnvironmentForLineno{#1}%
  \patchAmsMathEnvironmentForLineno{#1*}}%
\AtBeginDocument{%
\patchBothAmsMathEnvironmentsForLineno{equation}%
\patchBothAmsMathEnvironmentsForLineno{align}%
\patchBothAmsMathEnvironmentsForLineno{flalign}%
\patchBothAmsMathEnvironmentsForLineno{alignat}%
\patchBothAmsMathEnvironmentsForLineno{gather}%
\patchBothAmsMathEnvironmentsForLineno{multline}%
}
\newtheorem{thm}{Theorem}[section]

\newtheorem{prop}[thm]{Proposition}

\newtheorem{defn}[thm]{Definition}

\newtheorem{example}[thm]{Example}

\def\env{{\rm env}}

\newcommand{\cT}[1]{\mathcal{T}}%{C(S^1)^{(#1)}}
\newcommand{\FR}[1]{\mathcal{FR}}%{C(S^1)_{(#1)}}

\def\V{\mathcal{V}}
\def\VR{E_{V_4}}%{C(S^1)_{(#1)}}

\def\D{\mathbb{D}}

\def\R{\mathbb{R}}

\def\I{\mathbb{I}}
\def\Z{\mathbb{Z}}
\def\C{\mathbb{C}}
\def\CP{\mathbb{CP}}

\def\dar[#1]{\ar@<2pt>[#1]\ar@<-2pt>[#1]}

\newcommand{\nil}[1]{}

\makeatletter
\DeclareMathOperator{\exterior}{\@ifnextchar^\@exterior{\@exterior^{}}}
\def\@exterior^#1{\mathop{\bigwedge\nolimits^{\!#1}}}
\makeatother

\begin{document}

\title{K-theoretic invariants of four three-dimensional operator systems}
\author{Walter D. van Suijlekom}

\address{Institute for Mathematics, Astrophysics and Particle Physics, Radboud University Nijmegen, Heyendaalseweg 135, 6525 AJ Nijmegen, The Netherlands.}

\email{waltervs@math.ru.nl}
\date{30 September 2026}

\begin{abstract}
  We compute $K$-theoretic invariants for four three-dimensional operator systems to illustrate the general theory. We find that in some cases the invariants distinguish the operator system from its $C^*$-envelope, while in others they do not. We also explore extensions of the theory. 
\end{abstract}

\maketitle

%\setcounter{tocdepth}{1}
%\tableofcontents

\section{Introduction}
We illustrate our recently proposed extension of K-theory from $C^*$-algebras to operator systems \cite{Sui24,Sui25b} by means of four three-dimensional examples where they can be computed.  
A fully detailed exploration and determination of the invariants for the Toeplitz operator system, its dual given by the Fej\'er--Riesz operator system and graph operator systems is contained in \cite{Sui26b}. The selection of low-dimensional examples in the current paper is meant to concisely illustrate the breadth and variety of possibilities. 

In one of the cases we consider here, there are more elements in the $K$-invariants of the operator system in relation to those of the $C^*$-envelope; in other cases it is the other way around, and in yet another case they coincide. We also consider refinements of the general theory.

\subsection*{Acknowledgments}
%We thank Kristin Courtney, Maximilian Illmer, Sam Kim, Yuezhao Li, Tim Netzer, Lukas Obermeyer and Hermann Schulz--Baldes for fruitful discussions concerning the derivation of the $K_0$-invariants for the operator systems presented here.
For some parts of this manuscript use has been made of LLMs. Specifically, there have been dialogues of the author with ChatGPT-5.6 Sol (OpenAI), and Claude Opus 5 (Anthropic), assisting in the search for the variety of examples presented here. This was then further developed, rephrased and improved by the author. Open weight LLMs (QWen3.8) were used for proofreading. That being said, the responsibility for the content below lies entirely with the author.

\section{General theory}
We recall from \cite{Sui24,Sui25b} the notion of non-singular elements and hermitian forms in an operator system. For our purposes it is sufficient to consider {\em concrete} operator systems, which are defined to be $*$-closed and norm-closed vector spaces of bounded operators on a Hilbert space and containing the unit 1. In fact, we realize an operator system $E$ inside its $C^*$-envelope via the map $\imath_E : E \to C^*_\env(E)$. For more details, we refer {\em e.g.} to \cite{ER00,Pau02,Pis03}

%Even though in an operator system we cannot speak about invertible elements, we may use the pure and maximal ucp maps to introduce the following notion of nondegeneracy. 
\begin{defn}
%\label{defn:nondeg}
\label{defn:non-sing}
Let $(E,e)$ be a unital operator system. An element $x \in M_m(E)$is called {\em non-singular} if $\imath_E^{(m)}(x)$ is invertible as a matrix with entries in the $C^*$-envelope $C^*_\env(E)$. 

A {\em hermitian form} is a non-singular element $x \in M_m(E)$ which is self-adjoint. 
\end{defn}
We will write $H(E,m)$ for all hermitian forms in $M_m(E)$ and $G(E,m)$ for all non-singular elements in $M_m(E)$. Also, set $\V_0(E,m) := \pi_0(H(E,m))$ and $\V_1(E,m) := \pi_0(G(E,m))$, {\em i.e.}, the sets of path-components of hermitian forms and non-singular elements in $M_m(E)$, respectively.  
%We write $x \sim_m x'$ if $x,x'$ belong to the same path component in either $H(E,m)$ or $G(E,m)$.
All the sets $\V_p(E,m)$ ($p=0,1$) are invariants of operator systems in the following sense. 
\begin{prop}
  If $E$ and $F$ are completely order isomorphic then $\V_p(E,m) \cong \V_p(F,m)$ for all $m \geq 1$.
\end{prop}

\begin{example}
  \label{ex:V-C}
  Consider $E=M_n(\C)$. Then $\V_0(M_n(\C),m)$ is the set of path components of invertible hermitian $nm \times nm$ matrices. Any such matrix $x$ can be diagonalized with a unitary matrix, and since the unitary group $U(nm)$ is connected, there is a continuous path between $x$ and the corresponding diagonal matrix. In turn, this diagonal matrix is path-connected (in the space of invertible hermitian matrices) to the corresponding signature matrix, which is unique up to ordering.  In other words, $\V_0(M_n(\C),m)$ can be parametrized by the signature of %either the number ${\bf n}_+$ of positive, or the number ${\bf n}_-$ of negative eigenvalues of
  the hermitian forms%. We choose the latter, so-called {\em negative index of inertia}
  , yielding an isomorphism
  \begin{equation}
    \label{eq:V-C-n}
  \V_0(M_{n}(\C),m) \cong \{  -nm , -nm+2 ,\ldots, nm-2, nm \}
 \end{equation}
\end{example}
Concerning the invariant $\V_1$, we record the following from \cite{Sui25b}:
 \begin{prop}
      \label{prop:K1}
Suppose that $E$ is a unital operator system with finite-dimensional $C^*$-envelope. Then for any $m \geq 1$ the space $G(E,m)$ is path connected. Consequently, $\V_1(E,m) = \{ [1^{\oplus m}] \}$ for all $m \geq 1$.
  \end{prop}
%\proof
%Suppose that $E$ has a finite-dimensional $C^*$-envelope $M_d(\C)$ (the general case follows by taking direct sums). We claim that $\imath_E(E) \cap M_d(\C)^\times$ is contractible to the point $\imath_E(e) = \mathbb{I}_d$ (and the same applies to $M_n(E)$). Namely, for any $x  \in \imath_E(E) \cap M_d(\C)^\times$ choose $z \in S^1$ so that no eigenvalue of $x$ lies on the ray through the origin and $z$ (since $x$ has finitely many eigenvalues, this is always possible).  For any $t \in [0,1]$ we define  $\gamma(t) = t x + z(1-t) \mathbb I_d$. Then $\det \gamma(t) = 0$ iff $z(1-1/t)$ is an eigenvalue of $x$. By our choice of $z$ we find that $\det \gamma(t) \neq 0$ for all $t \in [0,1]$. Moreover, since complex linear combinations of $x$ and $ \mathbb I_d$ are contained in $\imath_E(E)$ it follows that the path $\gamma(t)$ is contained in $\imath_E(E) \cap M_d(\C)^\times$. More generally, we may conclude that $\imath_E^{(n)} (M_n(E)) \cap M_{nd}(\C)^\times$ is contractible to a point for any $n \geq 1$. 
 %\endproof
The limit structure $\varinjlim_m \V_p(E,m)$ was analyzed in \cite{Sui24,Sui25b} to come to a definition of $K$-theory for operator systems, which was furthermore shown to be invariant under Morita equivalence of operator systems. In the low-dimensional examples treated in this paper, it will be sufficient to restrict ourselves to analyzing $\V_p(E,m)$ for the case $m=1$. In order for the reader not to get lost in notation, we will therefore also write $G(E,1) \equiv E^\times$ and $H(E,1) \equiv E_h^\times$. Accordingly we have $\V_0(E,1) = \pi_0(E_h^\times)$ and $\V_1(E,1) = \pi_0(E^\times)$.

\section{Four three-dimensional operator systems}
\subsection{Operator system of Toeplitz matrices}
\label{sect:toep-pi1}
We consider the space $\cT{2}$ of $2 \times 2$ complex Toeplitz matrices, and write $T \in \cT{2}$ as
$$
T = \begin{pmatrix} t_0 & t_1 \\ t_{-1} & t_0 \end{pmatrix} ; \qquad (t_{-1},t_0,t_1 \in \C).
$$
The $C^*$-envelope of $\cT{2}$ is $M_2(\C)$ as can be easily seen. 

\begin{figure}
  \includegraphics[scale=.5]{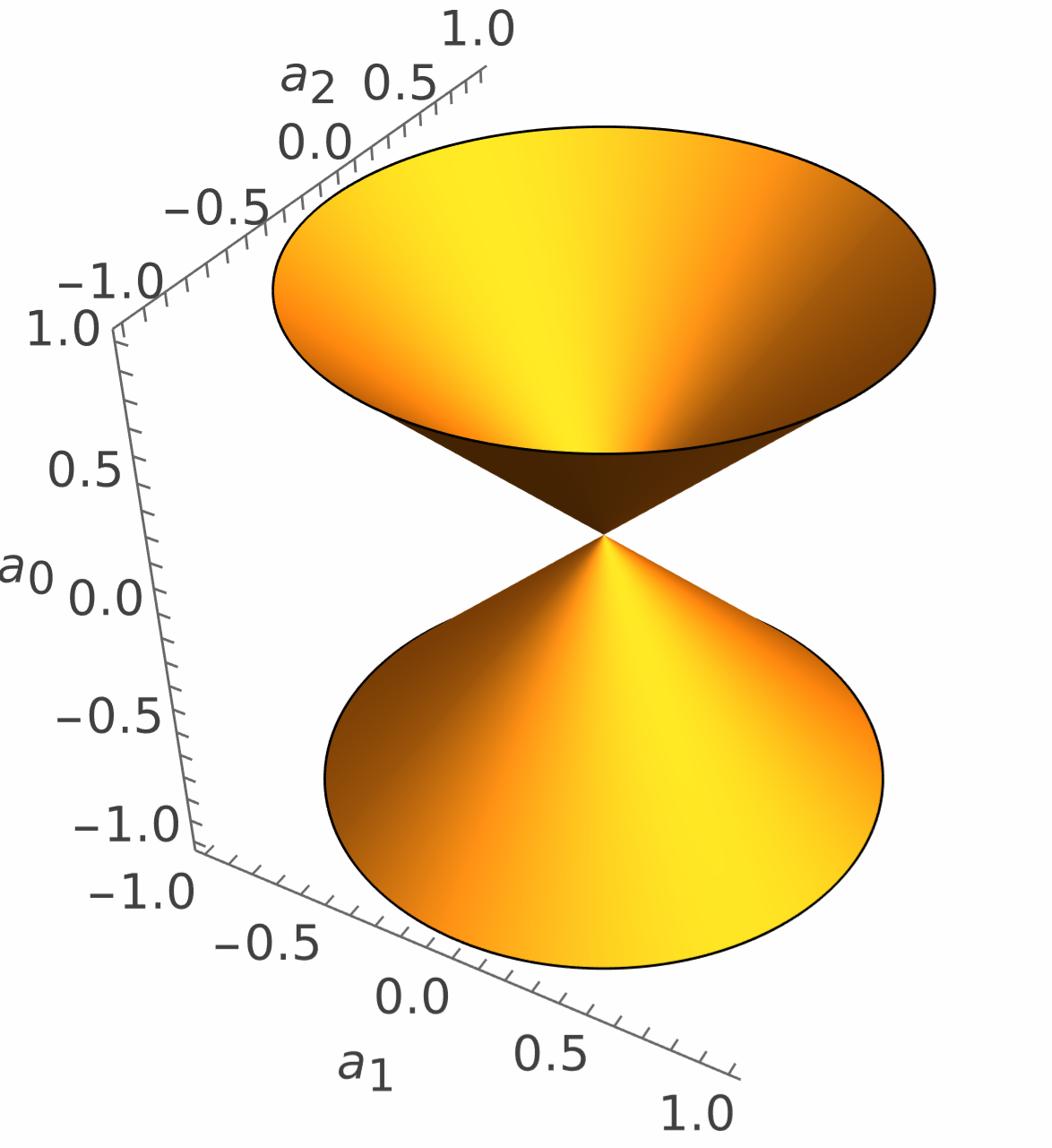}
  \caption{The singular locus of $a_0^2 -a_1^2 -a_2^2$ as a function of $(a_0,a_1,a_2) \in \R^3$; its complement is the space $\cT{2}_h^\times$ ({\em cf.} subsection \ref{sect:toep-pi1})}
  \label{fig:toep}
  \end{figure}

For the invariant $\pi_0(\cT{2}_h^\times)$ we should consider invertible hermitian Toeplitz matrices, say, 
$$
T = \begin{pmatrix} a_0 & a_1+i a_2 \\ a_1-i a_2 & a_0 \end{pmatrix}; \qquad (a_0,a_1,a_2 \in \R). 
$$
The set of invertible elements is given by the complement of the singular locus; the latter is the cone defined by $\det(T) = a_0^2-a_1^2 -a_2^2=0$ (see Figure \ref{fig:toep}).  
This cone divides $\R^3$ into three path-connected regions, which amounts to the fact that the path components of $H(\cT{2},1)$ are specified by matrix signature as explained above. In other words,
$$
\pi_0(\cT{2}_h^\times) \cong \{ -2,0,2\}.
$$
We already know from Example \ref{ex:V-C} that the same applies to the $C^*$-envelope, {\em i.e.}  $\pi_0(M_2(\C)_h^\times) \cong \{ -2,0,2\}$.

Concerning the invariant $\pi_0(\cT{2}^\times)$ we can be brief, as it follows from Proposition \ref{prop:K1} that it is trivial, and coincides with the invariant of its $C^*$-envelope.

As such the $K$-invariants cannot distinguish the $2 \times 2$ Toeplitz operator system from its $C^*$-envelope. Interestingly, however, one can also conclude from the above geometric description that
$$
\pi_1(\cT{2}_h^\times) = \Z , \qquad \pi_k(\cT{2}_h^\times) =0; \qquad (k > 1).
$$
In contrast, for the $C^*$-envelope $M_2(\C)$ the space $M_2(\C)_h^\times$ splits into two definite parts and one indefinite part (of matrix signature 0). The first are open convex cones, hence contractible, while the indefinite part by Sylvester's Law is a homogeneous space for $GL_2(\C)$ acting by conjugation on $\Lambda = \left(\begin{smallmatrix} 1 & 0 \\ 0 &  -1 \end{smallmatrix} \right)$. Thus, this part is homeomorphic to $GL_2(\C)/ U(1,1)$ which deformation retracts to $U(2)/U(1)\times U(1) = \CP^1$. We conclude that
$$
\pi_1( M_2(\C)_h^\times) = 0, \qquad \pi_2(M_2(\C)_h^\times ) = \Z, \qquad 
 \ldots
$$
In conclusion we can say that the homotopy groups of the space of hermitian forms do distinguish the Toeplitz operator system from its $C^*$-envelope; in fact, already the fundamental group does. This forms an interesting connection to the homotopy invariants of Banach algebras considered in {\em e.g.} \cite{CL86,Rie87}.

\subsection{Fourier truncations on the circle}
Next, we consider the Fej\'er--Riesz operator subsystem $\FR{2}$ of $C(S^1)$ (considered at length in \cite{CS20}) , which is spanned by the generator $z$, its adjoint $z^*$ and the unit. In other words, we consider functions on $S^1$ of the form
$$
a(z) = a_{-1} \overline z + a_0 + a_1 z ; \qquad (z \in S^1). 
$$
Hermitian forms in $\FR{2}$ are given by real-valued functions which are invertible on $S^1$, in particular, they are either strictly positive, or strictly negative. The corresponding two cones in the vector space $\FR{2}$ are of course contractible, so that $\pi_0(\FR{2}_h^\times) \cong \{\pm 1\}$. The same is actually true for $C(S^1)$ so that also $\pi_0(C(S^1)_h^\times) \cong \{ \pm 1 \}$. In \cite{Sui26b} we show that this continues to hold for the invariant $\V_0(\cT,m)$ for any $m$. 

Instead, we consider the invariant $\pi_0(\FR{2}^\times)$, based on the space of invertible elements in $\FR{2}$. This amounts to the functions $a(z)$ being nowhere vanishing (on $S^1)$. Let us write
$$
a(z) = \overline z  (a_{-1} + a_0 z + a_1 z^2) .
$$
Then the polynomial $(a_{-1} + a_0 z + a_1 z^2)$ has two zeros which lie inside or outside of $S^1$. Accordingly we write: for the possible cases:
$$
a(z) = \left\{ \begin{array}{ll}
  c \cdot \overline z (1-\gamma z) (1-\gamma' z) \sim \overline z\\
  c' \cdot \overline z (z-\alpha)(1-\gamma z) \sim 1 \\
  c'' \cdot \overline z (z-\alpha) (z-\alpha') \sim z 
 \end{array} \right .
$$
for $\alpha,\alpha',\gamma,\gamma' \in \D$, by moving these numbers to $0$. These possibilities correspond to winding numbers $-1,0$ and $1$, respectively. Since the winding number is a homotopy invariant, we conclude that $\pi_0(\FR{2}^\times) \cong \{-1,0,1\}$. This should be contrasted to $\pi_0(C(S^1)^\times) \cong \Z$, since in that case all winding numbers occur.

\subsection{An operator system related to the Vierergruppe}
\label{sect:more-env}
Let $V_4 = \Z_2 \times \Z_2$ and consider the unital operator system $\VR$ generated by the non-trivial characters $\chi_1, \chi_2$ given by $\chi_1(s,t) = s$ and $\chi_2(s,t) = t$. In other words, $\VR$ consists of functions in $C(V_4) \cong \C^4$ of the following form:
$$
f(s,t) = a_0 +  s a_1 +t a_2 ; \qquad (a_0,a_1, a_2 \in \C)
$$
One has $C^*_\env(\VR) \cong C(V_4) \cong \C^4$ because for each of the four points of $V_4$ one can find functions in $\VR$ that peak at that point, {\em e.g.} 
$$
f_{+,+}(s,t) =2 + s + t \leadsto  \| f_{+,+} \| = 4,
$$
which is peaking at $(1,1)$. Similarly for the other points so none of the summands in $C(V_4)$ is a boundary ideal for $\VR$.  

\begin{figure}
\begin{subfigure}[l]{0.4\textwidth}  
  \includegraphics[scale=.5]{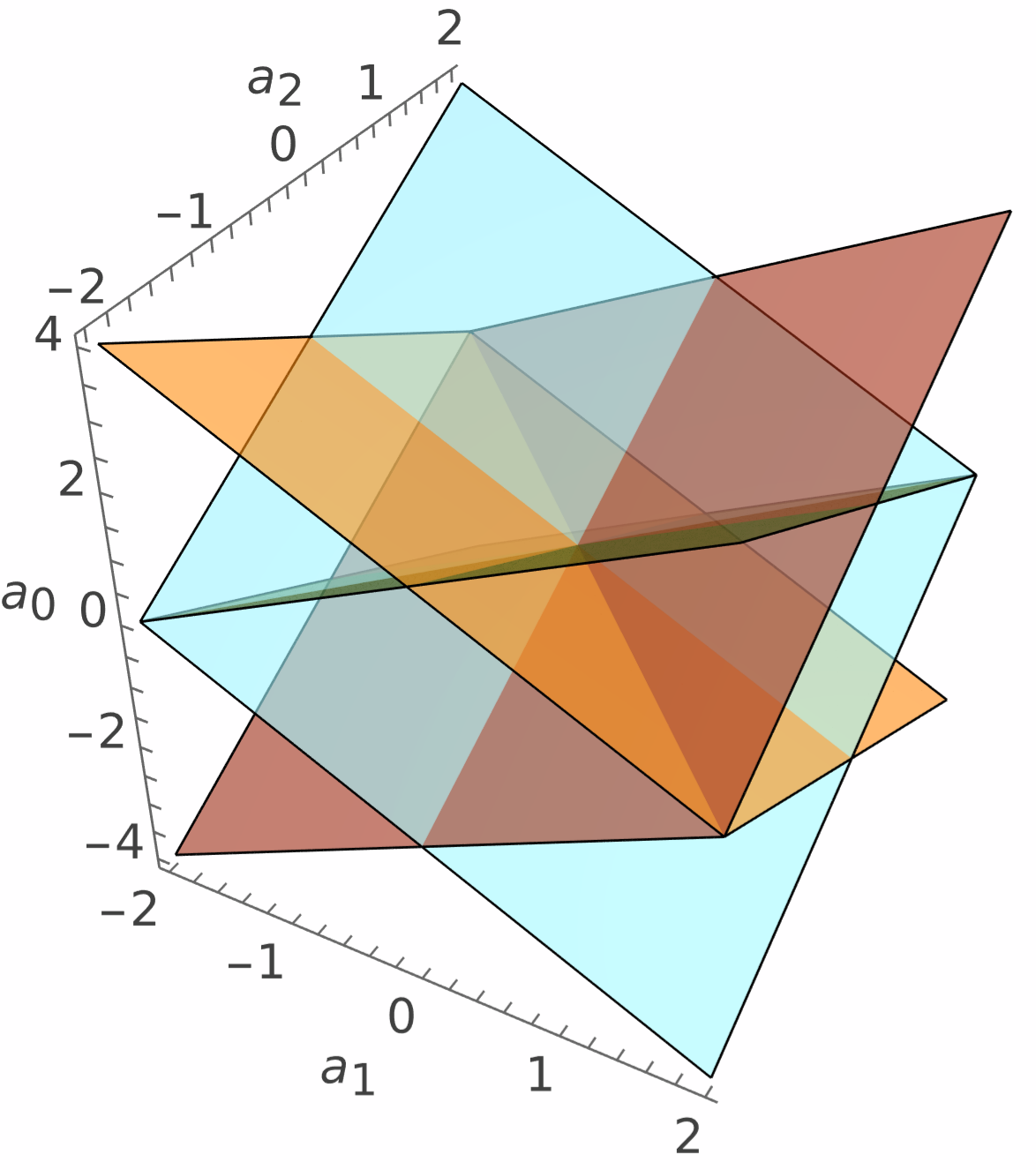}
  \caption{The singular locus of $a_0 + a_1 \chi_1  + a_2  \chi_2$ as a function of $(a_0,a_1,a_2) \in \R^3$; its complement is the space $(\VR)_h^\times$ with 14 path components ({\em cf.} subsection \ref{sect:more-env})}
  \label{fig:14conn}
\end{subfigure}\hspace{1cm}
\begin{subfigure}[r]{0.4\textwidth}  
  \includegraphics[scale=.5]{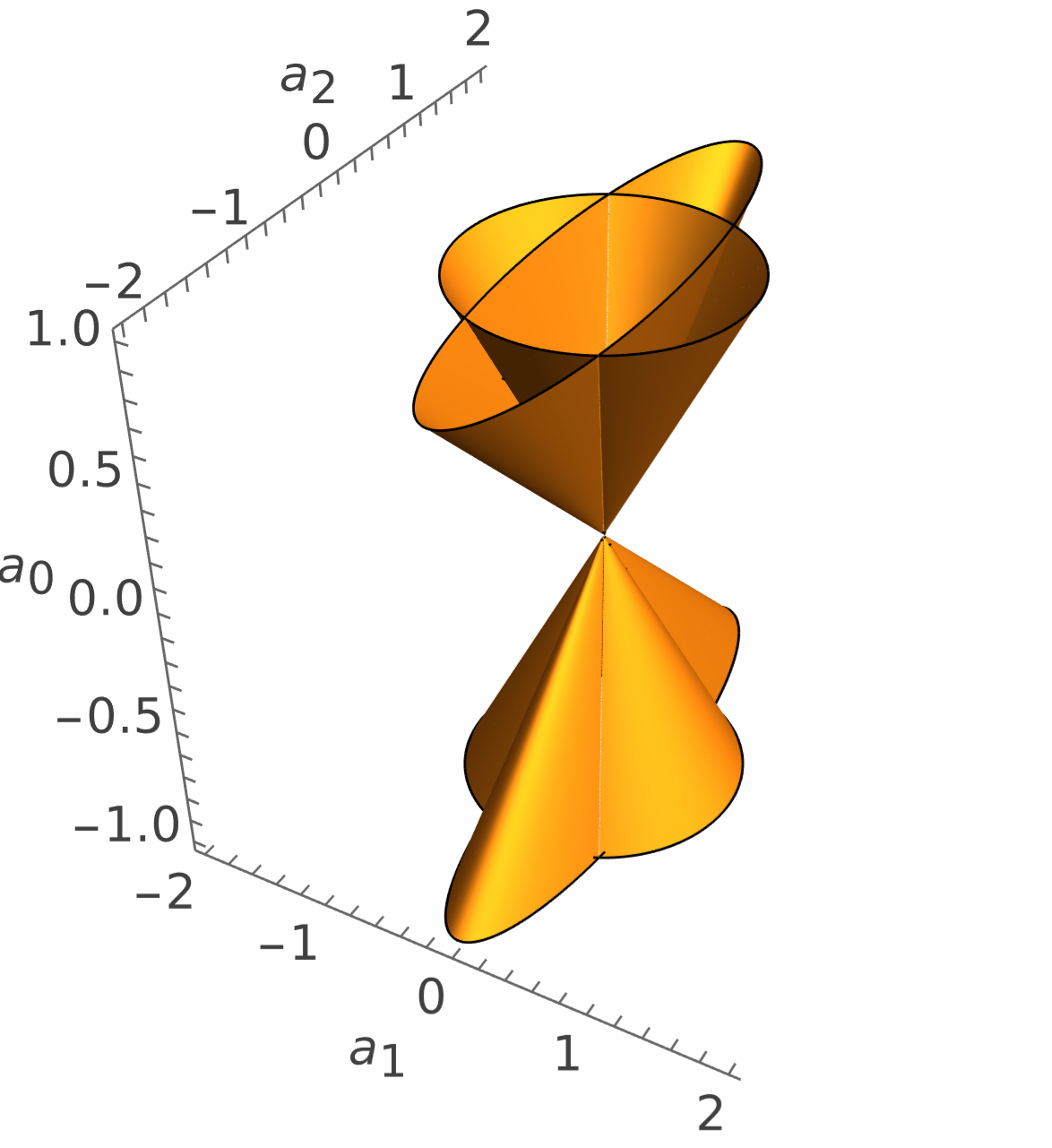}
  \caption{The singular locus of $a_0 \I_2 \oplus \I_2  + a_1 X + a_2 Y$ as a function of $(a_0,a_1,a_2) \in \R^3$; its complement is the space $(E_{\sigma})_h^\times $ with 11 path components ({\em cf.} subsection \ref{sect:more-os})}
  \label{fig:11conn}
\end{subfigure}
\caption{The two operator systems considered in Sections \ref{sect:more-env} and \ref{sect:more-os} for which the invariant $\V_0(E,1)$ can distinguish the operator system from its $C^*$-envelope.}
\end{figure}

The invariants $\pi_0(C(V_4)_h^\times)$ of the $C^*$-envelope can easily be determined to be
$$
\pi_0(C(V_4)_h^\times) = \pi_0 ((\R^\times)^4) \cong \{ \pm 1 \} ^4
$$
%corresponding to the 16 possible signs at the four points of $V_4$.
In contrast, an element in $(\VR)_h^\times$ is a nowhere-vanishing real-valued function $f(s,t)$ satisfying 
$$
f(1,1) + f(-1,-1) = f(1,-1)+f(-1,1).
$$
But this implies that a function $f$ for which both $f(1,1)$ and $f(-1,-1)$ are strictly positive, while $f(1,-1)$ and $f(-1,1)$ are strictly negative, cannot lie in $(\VR)_h^\times$. 
An illustration of this is given in Figure \ref{fig:14conn}. The signatures $(1,-1,1,-1)$ and $(-1,1,-1,1)$ at the four points of $V_4$ cannot be reached by a real linear combination of $1,\chi_1$ and $\chi_2$. In fact, the vanishing of $a_0 + s a_1 + t a_2$ divides Euclidean 3-space into 14 path components forming $(\VR)_h^\times$. Hence, $|\pi_0((\VR)_h^\times)|=14$, while $|\pi_0(C(V_4)_h^\times)| = 16$.

\subsection{An operator system spanned by Pauli matrices}
\label{sect:more-os}
We consider $A= M_2(\C) \oplus M_2(\C)$ and take the operator subsystem
$$
E_{\sigma}:= \text{span}_\C \{\I_2 \oplus \I_2, X , Y\};\qquad  X= \sigma_1 \oplus 2 \sigma_1, \qquad Y= \sigma_2 \oplus \frac 12 \sigma_2
$$
One checks that $A$ is the $C^*$-envelope: none of the blocks $M_2(\C)$ is a boundary ideal as they strictly reduce the norm of one of the summands in $E_\sigma$. Hence, $(E_\sigma)_h^\times$ consists of linear combinations $a_0 \I_2 \oplus \I_2  + a_1 X + a_2 Y$ with $(a_0,a_1,a_2) \in \R^3$ such that the resulting matrix is non-singular. One computes
$$
\det (a_0 \I_2 \oplus \I_2  + a_1 X + a_2 Y) = (a_0^2-a_1^2-a_2^2)(a_0^2-4 a_1^2 - \frac 14 a_2^2)
$$
so that $(E_\sigma)_h^\times$ is given by the complement of the singular locus.
This singular locus divides $\R^3$ into 11 different path-components: for a base of the cone $a_0=1$ or $a_0 =-1$ we find that it is formed by a circle and an ellipse in the coordinates $(a_1, a_2)$ (see Figure \ref{fig:11conn}). In contrast, $(M_2(\C) \oplus M_2(\C))_h^\times$ has 9 path-components ({\em cf.} Example \ref{ex:V-C}).

\section{Outlook}
We have considered the $K$-theoretic invariants $\V_0(E,1)$ and $\V_1(E,1)$ for four three-dimensional operator systems. 
In these examples we determined whether or not the invariant distinguishes the operator system from its $C^*$-envelope. We summarize the possible outcomes for these examples in the table in Figure \ref{fig:table}, leaving the refinement of the general theory by including higher homotopy groups to future work. 
A full investigation of the invariants $\V_p(E,m)$ for $p=0,1$ and any $m\geq 1$ for the operator systems of $n \times n$ Toeplitz matrices, for its dual given by the Fej\'er--Riesz operator system, and all graph operator systems is contained in \cite{Sui26b}.

\begin{figure}

  \begin{tabularx}{.9\textwidth}{|c||c|c||c|c||c|c||c|c|}
    \hline
    & $\cT{2}$ & $ M_2(\C)$ & $\FR{2}$ & $C(S^1)$ &$ \VR$ & $\C^4$ & $E_\sigma
    $& $M_2(\C) \oplus M_2(\C)$ \\
    \hline
    $\pi_0(E_h^\times)$ & \multicolumn{2}{c||}{ $\{ -2,0,2\}$ }& \multicolumn{2}{c||}{$\{ \pm 1\}$}&\cellcolor{lightgray} $\#: 14$ &\cellcolor{lightgray} $\#: 16$ & \cellcolor{lightgray}$ \#: 11 $ &\cellcolor{lightgray} $\#: 9$
    \\
    \hline
    $\pi_0(E^\times)$ &\multicolumn{2}{c||}{$ \{ 1 \}$}&\cellcolor{lightgray} $\{ -1,0,1 \} $ &\cellcolor{lightgray} $\Z$&\multicolumn{2}{c||}{$ \{ 1 \}$}&\multicolumn{2}{c|}{$ \{ 1 \}$} \\
    \hline \hline
    $\pi_1(E^\times_h)$ &\cellcolor{lightgray} $\Z$ &\cellcolor{lightgray} 0 \\
    
%    $\pi_1(E^\times)$ & $\Z$ & $\Z$ \\
    \cline{1-3}
    \end{tabularx}
  \caption{$K$-invariants for the three-dimensional operator systems considered in this paper; both $\pi_0(E_h^\times)$ and $\pi_0(E^\times)$ are indicated for both the operator system and its $C^*$-envelope. Gray indicates that the $K$-invariant distinguishes the operator system $E$ from its $C^*$-envelope.}
  \label{fig:table}
  \end{figure}

%\bibliographystyle{plainmath}
%\bibliography{references}

\newcommand{\noopsort}[1]{}

\end{document}